\documentclass{article}
\usepackage{graphicx} 
\usepackage{amsmath,amsfonts,amssymb,amsthm}
\usepackage{physics}
\usepackage{geometry}
\usepackage{hyperref}
\numberwithin{equation}{section}
\usepackage{indentfirst}
\newtheorem{definition}{Definition}
\newtheorem{theorem}{Theorem}
\newtheorem{lemma}{Lemma}
\usepackage{mathrsfs}

\title{A Stochastic Approach to the Definition of the Path Integral Measure}
\author{Timur Obolenskiy  \\\href{mailto:me@somewhere.com}{timurobolenskiy@gmail.com} }

\date{December 25, 2025}

\begin{document}

\maketitle

\section{Introduction}

The Path Integral was first introduced by Feynman in 1948, reminiscent of earlier works by Wiener that mainly sought to describe stochastic dynamics related to Brownian Motion. In its roots, the Path Integral stood as a heuristic construction, motivated physically as the summation over all possible trajectories. The Heuristic physical foundation for the Path Integral consists of summing over the contributions of all possible trajectories a particle may take. While this view may be useful from an intuitive standpoint, it is problematic, as the summation over all trajectories is left undefined and the measure $\mathcal{D}x(t)$ is impossible to ground in rigorous mathematics due to Cameron's theorem, which proves that the Feynman Measure is not a measure. As such, the root of the problem of formulating a rigorous path integral measure may be seen as being a problem of the path space being considered. 
\vspace{1em}

Over time, several attempts have been made at producing a rigorous formulation of the path integral, with the prominent Feynman–Kac theorem linking Euclidean path integrals to diffusion semigroups of stochastic processes, providing a rigorous framework that, however, is only applicable after Wick rotation and hence is set in imaginary time. Further refinements by Cameron, Martin and others have all made progress in placing the Euclidean theory on solid mathematical footing, yet remained extendable to real-time Quantum Mechanics only through the means of Analytic Continuation. In Lorentzian signature, the difficulties of oscillatory integrals, alongside the non-existence of the measure and infinite-dimensional geometry persisted.
\vspace{1em}

In the present work, we provide a possible approach to resolving the issues of convergence and dimensionality by restricting the relevant geometry to a neighborhood of the classical trajectory, which through fibration also provides an avenue for establishing a Gaussian measure for integration, yielding a controlled, well-defined functional integral over a stochastic process. Thus, the difficulties of the path integral stem not only from oscillatory factors but from the choice of path space itself, particularly the ill-defined “sum over all paths” which includes trajectories that are physically inaccessible to the system. Classical mechanics already suggests strong restrictions: for a fixed energy, the Euler–Lagrange equations determine a unique trajectory between two points on a configuration manifold under certain boundary conditions (see Lemma 1). Even when the energy varies within a small uncertainty window, admissible paths remain confined to a bounded tubular neighborhood of the classical trajectory. Extrapolating this to the quantum domain, the finite spread of a wave packet and the uncertainty relation imply that only trajectories within a controlled deviation from the classical path contribute significantly. Thus, a physically meaningful reformulation of the path integral may be built on a restricted path space $C([0,1],M)$ that encompasses all physically relevant (accessible) trajectories.
\vspace{1em}

In the present paper, we attempt to develop a rigorous mathematical framework on such a restricted path space by considering a quadratic Lagrangian for a non-degenerate system and reformulating the trajectories through a stochastic process on bundles defined on a tubular neighborhood of a classical trajectory. By constructing a tubular neighborhood 
and expressing deviations as sections of the normal bundle, we obtain a restricted path space on which integration proceeds in a well-defined manner. The view of the tubular neighborhood arises naturally when considering a free particle. If one considers the alternative paths taken by a particle while traveling from A to B and accounts for the boundedness in maximal allowed deviation from the classical trajectory, imposed by finiteness of energy, then we obtain a trajectory with a peak at a given point whose magnitude we may denote $\eta$. However, the peak may be displaced along the whole classical trajectory, and even rotated, yielding a cylinder with radius $\eta.$ Within this geometric setting, we define a stochastic path integral as a Gaussian expectation weighted by the physical action, and we show that it admits a semigroup structure analogous to the Feynman–Kac representation.
\vspace{1em}

Finally, for full physical and mathematical equivalence, we prove that the Stochastic Path integral holds under Wick transformation and, in fact, reproduces the results of  the Feynman-Kac Theorem, which allows us to carry over the vast array of known results and properties of Euclidean path integrals onto our stochastic framework. In the current paper, we will only view the case when 1) a classical trajectory is given with fixed endpoints, finite time and finite energy uncertainty 2) the particle surely starts and ends at the endpoints of the classical trajectory, 3) any trajectory not ending at the endpoint of the classical trajectory does not exist, 4) the particle is non-relativistic. These conditions are technicalities imposed for simplicity of formulation: nothing in the framework forbids extensions to the domains where those conditions fail, however that necessitates a lot more work and refinement for mathematical consistency and rigor. We will be analyzing the case of the Gaussian free particle in this paper, yet the stochastic structure leveraged for omitting trajectories explicitly in integration will be formulated as useful for the general case, accounting for drift and more. The evolution of the system accounting for energy uncertainty, degenerate trajectories and infinite timescales will not be explored.
\vspace{1em}

\section{Redefining the Path Space}
\vspace{1em}

Similarly to how gauge redundancy leads to divergences due to contributions of infinitely many physically irrelevant quantities, it is possible to imagine that the contributions of physically irrelevant trajectories also contributes to the mathematical hurdles of the Path Integral in the particle theory. The summation over all paths is not only impossible to formulate under some standard measure, but also makes us account for paths which are impossible for the trajectory to take.
\vspace{1em}

\subsection{The Classical View}
\vspace{1em}
Let us consider the case of a particle traveling from point A to point B on a smooth, connected configuration manifold $M$. Then, the action for any path $\gamma$ is given by the integral of the Lagrange function $L(q,\dot{q},t)$ over the time-interval $[t_1, t_2]$. Given that all paths can essentially be represented by curves on $M$, we may parametrize them to our convenience and as such will consider the time-interval to be $[0,1]$ such that $ \forall \gamma  |_M, 
\space \gamma:[t_1,t_2] \rightarrow M, \space \gamma(t)=\gamma(x), \space x \in [0,1], \space t_1=0, t_2=1$ and $\gamma(t_1)=\gamma(0)=A, \space \gamma(t_2)=\gamma(1)=B$. Then the Path Space $\mathcal{P}=C^1([0,1],M)$ admits a unique $\gamma$ for all energies of the system. This can be formalized in the following lemma;
\vspace{1em}

\begin{lemma}
$\forall E \in \sigma$ where $\sigma$ is the non-degenerate part of the energy spectrum for a physical system with a quadratic Lagrangian $L(q,\dot{q}), \quad \exists\gamma_E \in \mathcal{P}, \space \gamma_E:[0,1] \rightarrow M, \space \gamma_E(0)=A, \space \gamma(1)=B.$ with $\gamma_E$ unique given $[0,T]$ contains no conjugate points and $T<\infty$.  
\end{lemma}
\vspace{1em}

\textbf{Proof:} A Manifold $M$ that satisfies the criteria laid out initially is any Riemannian Manifold $(M,g)$. Let the coordinates of this manifold be $q=(q^1,...,q^n)$. The Lagrangian whose explicit time dependence is eliminated and the corresponding action functional are given by
\begin{equation}
    L(q,\dot{q})=\frac{1}{2}g_{ij}(q)\dot{q}^i\dot{q}^j-V(q) \qquad S(\gamma)=\int_0^1L(q(t),\dot{q}(t))dt
\end{equation}
We may now introduce the canonical momentum $p_i=\frac{\partial L}{\partial \dot{q}^i}=g_{ij}\dot{q}^j$ with the corresponding Hamiltonian being $H(q,p)=\frac{1}{2}g^{ij}(q)p_ip_j+V(q)$ which satisfies the HJE, 
\begin{equation}
    \frac{\partial S}{\partial t}(q,t)=-H(q, \nabla_qS(q,t)), \quad \nabla_qS=g_{ij}(q)\frac{\partial S}{\partial q^j}\frac{\partial S}{\partial q^i} \qquad S_E(q,t)=\Phi_E(q)-Et
\end{equation}
The value of energy $E$ is given by the inverse-map from phase space from $M$ of $H 
\in C^\infty $, while $\Phi_E:M \rightarrow \mathbb{R}$ satisfies $H(q, \nabla W_E(q))=E$. The trajectory $\gamma_E$ that corresponds to the path taken by the system given it having energy $E$ is then obtained by integrating $\dot{q}^i(t)=\frac{\partial H}{\partial p_i}(q, \nabla W_E(q)), \space q(t_0)=A.$ Given that $H$ is convex in $p$, the HJE determines a unique element in $C^1$ and hence ascribes a unique $\gamma $ to all energies. Since $E$ is non-degenerate, there isn't a configuration of the system with the same energy but different "state". For example, if we have two objects with the same energy but one has some non-zero angular momentum, that second object is excluded from our analysis. Since there are no conjugate points contained in the interval, $\delta^2S$ is at all times non-degenerate, alleviating the issue of local locus-straying (minimal deviations induced by the Jacobi field). Hence, uniqueness is guaranteed. \qed
\vspace{1em}

For the entirety of the paper, we will be assuming systems under the conditions of Lemma 1. The uniqueness of trajectories for any given fixed energy value is a rather trivial result when discussing the boundary conditions imposed in Lemma 1, yet it is valuable to gain perspective into the nature of classical dynamics for non-degenerate systems. Importantly, unless the boundary conditions are exactly those imposed for Lemma 1, the plurality of trajectories becomes possible. Problems such as the non-uniqueness of trajectories for $T\rightarrow \infty$ and degenerate energies will not be explored in this paper for that would necessitate a global analysis. Now, let us assume the case where the energy isn't strictly fixed, varying under some small uncertainty values. The uncertainty in energy will also be considered as finite, as the analysis of systems with infinite uncertainties again collapses to the case of studying the global evolution for the system over the entire energy spectrum, which is outside of the scope of the present paper. In our case, we may consider a Lie Algebra of transformations of some base trajectory $\gamma_0$ under parameters $\alpha$. We may evaluate this through symplectic geometry on $T^*M$ (Phase Space).
\vspace{1em}

Let $\theta=p_idq_i, \space \omega=dp_i \wedge dq^i: \space \exists X_H, \quad \iota_{X_H}\omega=dH$ define the Hamiltonian field. Then, $\exists \mathcal{J}^t: T^*M \rightarrow T^*M, \space (\mathcal{J}^t)^*\omega=\omega$ gives a Hamiltonian flow preserving our symplectic 2-form. Let $I_\Gamma$ be the Poincare Invariants for our system. Then, for a Lie Group $G_\alpha$ of symplectomorphisms performed by transformation operators $T_\alpha$ on $T^*M, \quad \Gamma^\alpha_E=T_\alpha\Gamma_E$, where $\Gamma \subset T^*M$ defines a closed contour. 
\vspace{1em}

The family $\{ \Gamma ^ \alpha _E \}: \alpha \in \mathbb{R} $ yields a Solenoid $\frak{B}$ of contours around the original trajectory with $\oint _{\Gamma^\alpha_E} \theta =\oint _{\Gamma_E} \theta$ a symplectic invariant. This solenoid is an invariant in $\mathcal{P}$, with all having an action $S(\gamma^\alpha_E)=\beta S(\gamma_0^
{\alpha_0}): \space \gamma_0^{\alpha_0}=T_\alpha^{-1}\gamma_0^
\alpha$. 

\subsection{Hamilton's Principle and Boundedness}
\vspace{1em}

The Hamilton Principle states that $\delta S=0$, wherefrom $S(q+\delta q, \dot{q}+\delta \dot{q})-S(q,\dot{q})=0$. Thus, for any trajectory inside $\frak{B}$, we may hypothesize the existence of some $C>0$ such that $|S(\gamma^\alpha)-S(\gamma_0)|\leq C$ for any trajectory that is sufficiently close yet not equal to $\gamma_0$ and some variance in energy of the system, which lies within $[E-\delta E, E+\delta E]$.

\begin{theorem}
    Let the path space be endowed with a Hilbert manifold structure such that $\mathcal{P}=H^1([0,1],M)$, and $||\gamma(t)||^2_{H^1}=\int_0^1(|q(t)|^2+|\dot{q}(t)|^2)dt$. Then, $\exists \mathcal{T}_r (\gamma_E)= \{ \gamma \in \mathcal{P} : ||\gamma - \gamma_E||_{H^1} \leq r \} $ be the tubular neighborhood of $\gamma_E$. Then, $\exists C >0: \Delta S =S(\gamma)-S(\gamma_E) , | \Delta S| \leq C \quad \forall \gamma \in \mathcal{T}_C$. 
\end{theorem}

The boundedness of variation for the action between different trajectories is standard for quadratic Lagrangians and may be verified by a local Taylor expansion. For any $C^2$ functional on a Hilbert Manifold, there exists a local bound in the neighborhood of a critical point, meaning that for the action $|S(x)-S(x_0)|\leq c||x-x_0||^2_H, \quad c\neq C$ with $||x-x_0|| \leq \varepsilon.$
\vspace{1em}

Physically speaking, any trajectory whose energy falls within the uncertainty range of the energy of some reference trajectory is at all times contained within some tubular neighborhood of the latter if we maintain that $\delta S=0$. If that was not true, then the transformations of coordinates under variation---expressible as transformations under $G_\alpha$ given in Lemma 1.---would not be  null and hence we would obtain a contradiction. The relevance of this theorem to the reformulation of path space lies in that it allows us to restrict the domain of integration to only those trajectories which can be physically taken. To fully establish a foundation for this in quantum theory, we will also look at how this boundedness translates into the energy uncertainty and wave-packet spread. 

\subsection{The Quantum-Mechanical View}
\vspace{1em}

In application to Quantum Mechanics, the principle of boundedness may be analyzed from the standpoint of the spread of the wave-packet that a particle corresponds to. As was first shown by Wu $\And$ Yang in \cite{PSI-BUNDLE} in application to Gauge theory, the wave function may be understood as sections of a fiber bundle. Specifically, given $\pi : M \times H \rightarrow M$ a complex Hilbert bundle, $\psi : M \rightarrow H, \space x \rightarrow \psi(x) $ represents a section. For probability densities $\rho (x) = \langle \psi (x), \psi(x) \rangle =|\psi(x)|^2 \in \mathbb{R}_{\ge 0}$ on $H$ with $\rho(\gamma(1))=1=\rho(\gamma(0))$ as boundary conditions (starting and ending at A,B) providing a fiberwise probability normalization $\forall \gamma \in\mathcal{T}_\eta$. We will leverage this view to analyze extension of Theorem 1. to the quantum domain. 
\vspace{1em}

Let $\hat{H}\psi =\sigma \psi $. We will be viewing a case of a 1-dimensional particle which may be very easily generalized to a multi-dimensional case. In the semiclassical view, the ansatz $\psi(x,t)=A(x,t)\exp(\frac{i}{\hbar}S(x,t)): -\frac{\partial S}{\partial t}=E(x,t)=\sigma(x,t)$ near the classical trajectory $\gamma_0.$ Let $T=t_2-t_1$ be the physical lifetime of $\psi$'s evolution from $A$ to $B$. Then, for the stationary case, $S_\sigma(x,t)=\Phi_\sigma(x)-\sigma t$, so along $\gamma$ with energy $E$, $S_E(\gamma)=S_\sigma(\gamma)+\sigma T-ET \implies \delta S(\gamma)=S_{\sigma+\delta E}(\gamma)-S_{\sigma}(\gamma)=-\delta ET$, which corresponds to the change in action due to change in energy from $\sigma $ to $\sigma+\delta E$. Given the finiteness of $\delta ET,$ we already get the hint that some bound of form of $C$ exists.
\vspace{1em}

The wave function for the Schrodinger equation with a time-dependent Hamiltonian is $\psi(x,t)=\exp(\frac{-i}{\hbar}\hat{H}t)\psi_0$ which upon Fourier transform to energy space yields an energy-space wave function $\phi(E)$ with 
\begin{equation}
    (\delta E)^2=\int (E-\sigma)^2|\phi(E)|^2dE, \qquad \phi(E)=\int\psi(x)\exp(\frac{i}{\hbar}Et)dx
\end{equation}
It thus follows that $|\delta S(\gamma)|=|S_{\sigma+\delta E}(\gamma)-S_\sigma(\gamma)|=|\delta E|T\leq C \implies |\delta E| = \frac{C}{T}$. Thus, $|E-\sigma|\leq \frac{C}{T}$. Heuristically, by the uncertainty principle, we get $\frac{1}{2}\hbar \leq C$. Here, like in semiclassical analysis, the diameter of our tube will be on the order of $\sqrt {2\hbar}$
\vspace{1em}

The role of the trajectory in quantum mechanics is often misconstrued due to a variety of metaphysical arguments. We are not going to engage in such discourse, as it is irrelevant to our discussion if the physical particle travels over all trajectories at once or singles one out as being preferable. It is also irrelevant to make any additional ontological statement beyond what is standardly accepted. Rather, it matters solely that the domain of all possible trajectories that may be taken given some specific energy---with the associated bounded variance---is bounded. The paths may be such that the particle takes them all at once, or perhaps only some bizarre combination. However what is important for us is that all possible trajectories that may physically be traveled remain bounded and are constrained from straying too far away from the classical trajectory. In other words, it is only the confinement of trajectories that is relevant to the path space, not the mode in which the particle takes them. Just as heuristically we may assert that the 'velocity' $\dot{\gamma}(t) \in [\dot{\gamma}(t)-\epsilon, \dot{\gamma}(t)+\epsilon]: \epsilon =O(\sqrt{\delta E}),$ the paths are such that $\exists c \in \mathbb{R}: c||\gamma-\gamma_0||_{H^1}^2 \leq ||S(\gamma)-S(\gamma_0)||_{H^1}\leq C\implies ||\gamma-\gamma_0||_{H^1} \leq \sqrt{\frac{C}{c}}$. On path space, all trajectories contributing to $K(B,A)$ respect the relation $||\gamma-\gamma_0||_{H^1} \leq \sqrt{\frac{1}{2}\hbar}$.

\section{Stochastic Processes on Path Space and Fibration}

\subsection{Geometry of the $\eta$-Tube and Bundle structure}
\vspace{1em}

Let $(M,g)$ be a compact smooth Riemannian manifold of class $C^\infty$. In our analysis, we will strictly view the case of a system detailed above, traveling from two fixed points, with stochasticity in between. There exists a classical trajectory, per Lemma 1, such that in $\mathcal{T}_\eta$ of that trajectory there are some alternate trajectories that gain relevance at the quantum level. This classical trajectory will be denoted as $\gamma_0  : [0,1] \rightarrow M \in C^2, \space \gamma_0(0)=A, \space \gamma_0(1)=B$.  Let $N \gamma_0 \rightarrow [0,1]$ be the normal bundle to $M$ along $\gamma_0$ defined as $ N_{\gamma_0(t)}\gamma_0 := \{ v \in T_{\gamma_0 (t)}M : g( \dot \gamma_0(t), v)=0 \}$
\vspace{1em}

The admissible perturbations (smooth homotopy class with pinned endpoints) are bounded to a finite neighborhood of the classical trajectory, meaning that on $M$ we may define a local normal coordinate chart, while the neighborhood is defined by an injectivity radius $\inf_{t \in [0,1]}\operatorname{inj}_M \gamma_0(t)=\eta$. This allows us to define a tube around $\gamma_0$ on M as 

\begin{equation}
    C=\{ \exp _{\gamma_0(t)}(v):v \in N \gamma_0, \space t \in [0,1], ||v||_g \le \eta \} \subseteq \mathcal{T}_\eta
\end{equation}

In the equation above, $\exp$ is used to denote a Riemann exponential, while the chart induced for the neighborhood of the classical trajectory will be defined as being Fermi-Walker in order to compensate for the potential curvature of $\gamma _0$. This tube may be defined in terms of a cylinder by taking a base $\sigma_1=B(A,\eta) \cap C, \space \sigma_2=B(B, \eta ) \cap C$. Let $H^1$ denote a Sobolev space. Then, we may define transverse perturbations ${\frak{p}}(t) \in H^1([0,1], N \gamma_0): \space \sup || {\frak{p}}(t) || \leq \eta $. Parametrizing $\frak{p}$ by $t$, we get ${\frak{p}}(0)={\frak{p}}(1)=0$. These perturbations may be grouped into a set $P_\eta$ to which we may associate the exponential chart $\Phi(t)$ defined as follows;

\begin{equation}
    \Phi:P_{\eta} \rightarrow  \Omega [A,B]:=\{ \gamma | H^1: \gamma(0)=A, \gamma(1)=B, \gamma \in C \}
\end{equation}

The properties of $P_\eta$ and $\Omega[A,B]$ include being complete, and for $\Omega$, also normed vector spaces. This may be formalized;

\begin{lemma}
Let $(M,g)$ be a Riemannian manifold with a fixed continuous reference path $\gamma_0\in C([0,1],M)$ defined over it. For $\eta>0$ and a supremum metric $\rho(x,y)=\sup d_g(x,y)$, $(P_\eta,\rho)$ is complete.
\end{lemma}

\textbf{\proofname:} The metric space $(C([0,1],M),\rho)$ is complete whenever $(M,d_g)$ is complete. Let us parametrize $\rho$ by some $t\in [0,1]$, with $\rho(x,y)=\sup_{t\in[0,1]} d_g(x(t),y(t))$. If $\gamma_k$ is $\rho$--Cauchy, then $\forall t\in [0,1]$, the sequence $\gamma_k(t)$ is $d_g$--Cauchy and hence converges in $M$; the pointwise limit $\gamma$ is continuous and $\gamma_k\to\gamma$ uniformly. Now, let $\gamma_k\subset P_\eta$ be $\rho$--Cauchy. By completeness of $C([0,1],M)$, $\gamma_k\to\gamma$ in $\rho$ for some $\gamma\in C([0,1],M)$. Since the map $\gamma \rightarrow \rho(\gamma,\gamma_0)$ is $\rho$--continuous, $\rho(\gamma,\gamma_0)=\lim_{k\to\infty}\rho(\gamma_k,\gamma_0)\le \eta,$ so $\gamma$ satisfies the tube constraint. By fixing $k$ and choosing some $r_k>0:
r_k<\text{inj}_M(\gamma_k(t))$ for all $t\in[0,1]$. For all sufficiently large $j$ we have $\rho(\gamma_j,\gamma_k)<r_k$, hence for each $t$ there is a unique minimizing geodesic from $\gamma_k(t)$ to $\gamma_j(t)$ depending continuously on $t$ on M. This defines a homotopy $H_{k,j}:[0,1]\times[0,1]\to M$ under fixed endpoints from $\gamma_k$ to $\gamma_j$. Passing to the limit $j\to\infty$ yields a homotopy $H_k$ from $\gamma_k$ to $\gamma$. Since $\gamma_k\sim\gamma_0$ relative to endpoints by transitivity, and hence $\gamma\in P_\eta$. Therefore, every $\rho$--Cauchy sequence in $P_\eta$ converges in $(P_\eta,\rho)$, i.e.\ $P_\eta$ is complete.
\vspace{1em}

\begin{lemma}
On the manifold M, we may define a Sobolev space pinned by A,B as $\Omega[A,B]:=\{\gamma\in H^{1}([0,1],M):\ \gamma(0)=A,\ \gamma(1)=B\}.$ $\Omega[A,B]$ carries a natural structure of a $C^\infty$ Banach manifold modeled on $H^{1}_{0}([0,1],\mathbb{R}^{n}):$  $n=\dim M$. If $\gamma_0\in\Omega[A,B]$ and $\eta>0$
is such that $\exp_{\gamma_0(t)}$ is a diffeomorphism on $B_\eta(0)\subset T_{\gamma_0(t)}M$ for all $t$,
then the tubular chart $\Phi:U_\eta\subset H^{1}_{0}([0,1],\gamma_0^{*}TM)\to \Omega[A,B],\qquad
\Phi(v)(t):=\exp_{\gamma_0(t)}(v(t)),$ is $C^\infty$ with inverse $\ell:\Phi(U_\eta)\to H^{1}_{0}([0,1],\gamma_0^{*}TM),\qquad
\ell(\gamma)(t):=\exp_{\gamma_0(t)}^{-1}(\gamma(t)).$

\end{lemma}

\textbf{\proofname:} We will again define everything relative to some $\gamma_0\in\Omega[A,B]$. By continuity of $t \rightarrow \gamma_0(t)$ and positivity of the injectivity radius along compact sets (which our intervals are, for example), $\exists \eta>0$ such that for every $t\in[0,1]$ the exponential map $\exp_{\gamma_0(t)}:B_\eta(0)\subset T_{\gamma_0(t)}M\to M$ is a $C^\infty$ diffeomorphism. Let $\gamma_0^{*}TM\to[0,1]$ be the pullback bundle and set $U_\eta:=\{v\in H^{1}_{0}([0,1],\gamma_0^{*}TM):\ \|v\|_{L^\infty}<\eta\}.$ Since $H^{1}([0,1],\mathbb{R}^{n})\hookrightarrow C^{0}([0,1],\mathbb{R}^{n})$, the $L^\infty$-constraint makes $U_\eta$ open in $H^{1}_{0}([0,1],\gamma_0^{*}TM)$. Then, let $\Phi$ be defined such that $\Phi(v)(t)=\exp_{\gamma_0(t)}(v(t))$. Then $\Phi(v)\in H^{1}([0,1],M)$ and $\Phi(v)(0)=A$, $\Phi(v)(1)=B$ because $v(0)=v(1)=0$, and hence $\Phi:U_\eta\to\Omega[A,B]$ is well-defined. In any local trivialization of $\gamma_0^{*}TM$ over $[0,1]$, the map $(t,w) \rightarrow \exp_{\gamma_0(t)}(w)$
is $C^\infty$ in $(t,w)$ on $[0,1]\times B_\eta(0)$, which by Nemytskii theorems in 1D
Sobolev spaces (composition theorems) \cite{NEM}, the induced map $v\mapsto \Phi(v)$ is $C^\infty$ as a map between Banach spaces $H^{1}_{0}\to H^{1}$. Since $\exp_{\gamma_0(t)}$ is a diffeomorphism on $B_\eta(0)$, the pointwise inverse
$\exp_{\gamma_0(t)}^{-1}$ yields $\ell(\Phi(v))=v$ and $\Phi(\ell(\gamma))=\gamma$ on $\Phi(U_\eta)$, and the same composition theorem gives that $\ell$ is $C^\infty$. For a general $\gamma\in\Omega[A,B]$, we may repeat the construction with $\gamma$ in place of $\gamma_0$ to obtain a family of charts modeled on $H^{1}_{0}([0,1],\gamma^{*}TM)$. On overlaps of the charts, transition maps are of the form $v \rightarrow \exp_{\gamma(t)}^{-1}\!\bigl(\exp_{\gamma'(t)}(v(t))\bigr)$ and hence are $C^\infty$ by the same Sobolev composition argument. Therefore, $\Omega[A,B]$ is a $C^\infty$ Banach
manifold, and in the tubular neighborhood of $\gamma_0$ the chart is $\Phi$ with inverse $\ell$. \qed
\vspace{1em}

Hence, the Riemann exponential gains an expression through this local chart, and will be denoted through the convolution of the perturbation as $\Phi ( {\frak{p}})(t) = \exp _{\gamma_0(t)}({\frak{p}}(t))$. From Lemma 1, we may define the Riemann logarithm  through the notation $\ell (\gamma)(t):= \log^R_{\gamma_0 (t)}{\gamma (t)} \in N_{\gamma_0(t)}\gamma_0$, defined on $H^1$ and pinned by $[A,B]$.
\vspace{1em}

Let $M$ now be also endowed with Levi–Civita connections $\nabla$, and let $\gamma:[0,T]\to M$ be a $C^1$ curve. Its pullback tangent bundle $\gamma^*TM
:=\bigl\{(t,v)\,\big|\,t\in[0,T],\,v\in T_{\gamma(t)}M\bigr\}
\xrightarrow{\ \pi\ }[0,T]$, with T denoting the physical lifetime of the particle in the sense of the time it takes for it to go from A to B along one of the paths. The Bundle over an $n$-dimensional manifold is a rank–$n$ real vector bundle whose fiber over $t$ is $T_{\gamma(t)}M$. The Riemannian metric $g$ and connection $\nabla$ pull back fiberwise to $\gamma^*TM$.
\vspace{1em}

At all points of the trajectory $\gamma(t) \in \gamma$, the particle may take a variety of different directions given all of them over the whole time span of traveling along any one of the trajectories end at B. The inherent stochasticity of dynamics in $\mathcal{T}_\eta$ may essentially be interpreted as being attributed to the ability of the system to deviate from any one path, choosing instead the direction of some vector, which yields another trajectory and hence allows for the repetition of a similar deviation. These deviations will be counted as long as after some $t_M\in [0,1]$, the system will only begin to take vectors which in sum converge to a trajectory $ \{\tilde{\gamma}:[0,1] \rightarrow M\} \subset \mathcal{T}_\eta:\tilde{\gamma}(1)=B.$ This may be modeled through vector fluxes at all points $\gamma(t)$, which through the definition of an inner-product with the metric of $M$, may be turned into a flux space.
\vspace{1em}

\begin{definition}[Flux space]
The \emph{flux space} along $\gamma$ is the Bochner $L^2$–space $F$ defined in (3.3), with $u(t)\in T_{\gamma(t)}M$ and $\lvert u(t)\rvert_g^2:=g_{\gamma(t)}\bigl(u(t),u(t)\bigr)$. The inner product on $F$ is $\langle u,v\rangle_{L^2}
:=\int_0^1 g_{\gamma(t)}\bigl(u(t),v(t)\bigr) dt.$ Elements of $F$ are fluxes, i.e. square–integrable vector fields along $\gamma$ for measurable $u: [0,1] \rightarrow \gamma^*TM$ and 
\begin{equation}
  F:=L^2\bigl([0,1];\gamma^*TM\bigr)=\Bigl\{u:\ \int_0^1 \lvert u(t)\rvert_g^2\,dt<\infty\Bigr\} , 
\end{equation}

\end{definition}

Each flux space may be endowed with a Hilbert structure. This is trivially true from the definition. The entire fiber bundle is then also endowed with a  Hilbert space structure, given that $
\forall t \in [0,1], \space \dim (T_{\gamma(t)}M) =n$, and the fibers being Hilbert spaces with the norm induced by $g$. The family $\{T_{\gamma(t)}M\}_{t\in[0,T]}$ is a measurable field of Hilbert spaces, and $F$ is the associated Bochner $L^2$–space. Completeness of $F$ then follows from completeness of each flux space over finite–dimensional fibers, yielding an $L^2$ isomorphic space over $\mathbb{R}$. For certainty sake, and precise isomorphism to $L^2$, we ought to also establish separability. 
\vspace{1em}

\begin{lemma}
Each flux space $F$ is separable.
\end{lemma}
    
\textbf{Proof}
Let there be a $C^1$ curve $\gamma:[0,1]\to M$ with an associated parallel orthonormal frame along $\gamma$ via an isometry $I:T _{\gamma(t_0)}M\to\mathbb{R}^n$ and parallel transport $P_{t_0\to t}$, $E_i(t):=P_{t_0\to t}I^{-1}e_i$.
 $\forall u\in F=L^2([0,T];\gamma^*TM)$, there is a unique expansion $u(t)=\sum_i u^i(t)E_i(t)$ with $u^i\in L^2([0,T];\mathbb{R})$, and the map $\Phi:F\to L^2([0,T];\mathbb{R}^n)$, $(\Phi u)(t)=(u^1(t),\dots,u^n(t))$ a linear isometry, since $|u(t)|*g^2=\sum_i|u^i(t)|^2$ and hence $|u|*{L^2}=|\Phi u|_{L^2}$.
Given $f\in L^2([0,T];\mathbb{R}^n)$, defining $u_f(t):=\sum_i f^i(t)E_i(t)$ yields $u_f\in F$ with $\Phi u_f=f$, meaning that $\Phi$ is an isometric isomorphism $F\cong L^2([0,T];\mathbb{R}^n)$.
The latter space is separable because the set of $\mathbb{Q}^n$–valued step functions with rational breakpoints is countable and dense in $L^2([0,T];\mathbb{R}^n)$. \qed
\vspace{1em}

\begin{definition}[Covariant Sobolev flux space]
Let the space $H^1(\gamma^*TM)
:=\Bigl\{u\in L^2([0,T];\gamma^*TM)\ \Big|\ 
\nabla_t u:=\nabla_{\dot\gamma(t)}u(t)\in L^2\bigl([0,T];\gamma^*TM\bigr)\Bigr\}$ be endowed with the norm $\|u\|_{H^1}^2
:=\|u\|_{L^2}^2+\|\nabla_t u\|_{L^2}^2.$ Then, $H^1(\gamma^*TM)$ is a Sobolev space isomorphic to all $F$.
\end{definition}
\vspace{1em}

In the parallel–transport frame $(E_i(t))$ one has the covariant derivative of all $u(t)$ be given by 
\begin{equation}
    \nabla_t u(t)
=\sum_{i=1}^n \dot u^i(t)\,E_i(t) \qquad \nabla_t E_i(t)=0
\end{equation}
It then follows that the norm of the $u$ on the Sobolev flux space is
\begin{equation}
    \|u\|_{H^1}^2
= \sum_{i=1}^n\bigl(\|u^i\|_{L^2}^2+\|\dot u^i\|_{L^2}^2\bigr)
\end{equation}
so the isometry $\Phi$ extends to an isometric isomorphism
\[
\Phi:H^1(\gamma^*TM)\xrightarrow{\ \cong\ } H^1\bigl([0,T];\mathbb{R}^n\bigr).
\]

It ought to be noted that if $(F_i(t))_{i=1}^n$ is any smooth moving frame along $\gamma$ with uniformly bounded connection coefficients such that $\nabla_t u = \dot u + \Gamma(t)u$ for some $\Gamma\in L^\infty\bigl([0,T];\mathrm{End}(\mathbb{R}^n)\bigr)$ in the corresponding coordinates, $\|u\|_{L^2}^2+\|\dot u\|_{L^2}^2$ on the coordinate functions is equivalent to the covariant $H^1$–norm
$\|u\|_{L^2}^2+\|\nabla_t u\|_{L^2}^2$. Concretely, there exist constants $0<c\le C<\infty$ with 
\begin{equation}
    c\bigl(\|u\|_{L^2}^2+\|\dot u\|_{L^2}^2\bigr)
\;\le\;
\|u\|_{L^2}^2+\|\nabla_t u\|_{L^2}^2
\;\le\;
C\bigl(\|u\|_{L^2}^2+\|\dot u\|_{L^2}^2\bigr).
\end{equation}
This is the one–dimensional vector–bundle analogue of norm equivalence in standard Sobolev spaces.
\vspace{1em}

\subsection{Stochastic Dynamics in $\mathcal{T}_\eta$}

In the present section and onwards, we will concern ourself with formulating the stochastic dynamics through which the various trajectories the particle's evolution occurs over may be defined without an explicit summation over them. We will formulate it in a general way that for a free particle may be seen as excessive. This section may be skipped entirely, as section 4 and 6 concern themselves purely with that which may be used for the free particle case. Naturally, for systems of a different sort--harmonic oscillators, particles in potential fields, etc--cannot be defined in a trivial tubular way that a free particle can, as the path space corresponding to those systems requires a different sort of compactification. Nevertheless, what applies for the stochastic dynamics determining trajectories here may be carried over to the general case, although possibly with some modifications.
\vspace{1em}

We may consider now a stochastic process in each flux space for modeling the random perturbations of homotopies of the classical trajectory. Through the stochastic dynamics in each fiber, we may obtain all possible paths that the particle may take. Let $\mathcal{T}_\eta(\gamma)\subset M$ denote the $\eta$–tube around $\gamma$, realized (for $\eta>0$ sufficiently small) as the image under the exponential map of the open disk bundle in the normal bundle $N\gamma\subset\gamma^*TM$, with the tube being given by 
\begin{equation}
\mathcal{T}_\eta(\gamma)
=\Bigl\{\exp_{\gamma(t)}(\nu)\,\big|\ t\in[0,T],\ \nu\in N_{\gamma(t)}\gamma,\ \lvert\nu\rvert_g<\eta\Bigr\}.
\end{equation}
Any trajectory $\tilde{\gamma}:[0,T]\to\mathcal{T}_\eta(\gamma)$ admits a unique representation $\alpha(t)=\exp_{\gamma(t)}(X_t),
\qquad  X_t\in N_{\gamma(t)}\gamma$, where $X\in F$ is a flux field along $\gamma$ taking values in the normal bundle $N\gamma$. Thus the space of all physically relevant perturbations is naturally modeled by a closed subspace of the flux space $F$. In this setting, a covariant stochastic dynamics of fluctuations around $\gamma$ is given by a $\gamma^*TM$–valued diffusion SDE
\begin{equation}
    dX_t=b_t(X_t)\,dt+\sigma_t \circ dW^\gamma_t 
\end{equation}
where $b_t$ and $\sigma_t$ act fiberwise on $T_{\gamma(t)}M$, $W_t^\gamma$ is a covariant Wiener process transported along $\gamma$, and the tube constraint $\lvert X_t\rvert_g<\eta$ enforces that $\alpha(t)=\exp_{\gamma(t)}(X_t)$ remains in $\mathcal{T}_\eta(\gamma)$. We may introduce a noise function $\xi (t)$ such that $\forall v \in  N_\gamma M$, we have $\sigma^2\nabla_t^{\,2}v=\xi(t)v+V_{E_c}(v)+\sigma W_t^\gamma$ with the 'energy cost', $E^v_c=\bigl|E(\gamma(t))-E(\exp_{\gamma(t)}(v))|$
, where $E(\exp_{\gamma(t)}(v))$ is the energy of displacement from $\gamma(t)$ along $v\in F_\gamma.$. Therefore, via trivialization and some $\Omega_{ij}(t)=\langle \bar{\nabla}_t e_i(t), e_j(t) \rangle$, where $\bar{\nabla}$– fibre gradient on the metric of the normal bundle, we obtain 
\begin{equation}
v(t)=\sum^{n}_{i=1}E_i(t)\chi_t^i 
\ \implies\ 
dv(t)=\sum^n_{i=1} E_i(t)\,d\chi _t^i+\sum^n_{i=1}\bar{\nabla}_t E_i(t)\,\chi_t^i\,dt
=\sum^n_{i=1} \Bigl(d\chi_t^i+\sum_{j=1}^n\Omega_{ij}(t)\chi^j_t\,dt\Bigr)E_i(t).
\end{equation}

\begin{equation}
\therefore\ 
d \chi_t^i=-\frac{1}{\sigma^2}\bigl(\chi^i(t)\xi(t)+V^\chi_{E_c}\xi(t)\bigr)\,dt
+\sigma\, dW^{\gamma,i}_t
-\Omega^i(t)\chi^i(t)\,dt,
\qquad
V^\chi_{E_c}=\frac{1}{\|\chi\|}\frac{\partial E_c}{\partial \chi }.
\end{equation}

The covariant field for the Wiener process may be defined as $B^C_t \in T_{\gamma(t)}M : \forall Y \in B^C_t,$ with
\begin{equation}
    dY_t=-\frac{1}{\sigma^2}\xi(t)Y_t\,dt-V_{E_c}(Y_t)\,dt+\sigma\, dB^C_t,
    \qquad
    dv_t=-\Bigl(\frac{1}{\sigma^2}\xi(t)v_t+\frac{1}{\sigma^2}V^v_{E_c}(t)\Bigr)\,dt+\sigma\, dW^\gamma_t.
\end{equation}

$B$ and $W$ are linked by the equation $B^C_t=\sum^n_{i=1}e_i(t)W^{\gamma, i}_t$ with the same noise $\xi(t)$.
\vspace{1em}

Alongside the equations above, we may formulate a general governing equation that incorporates the necessary stochasticity that is analogous to the Euler-Lagrange equation, given by 
\begin{equation}
    \sigma^2\nabla_t^{\,2}Y_t=\xi(t)Y_t+V_{E_c}(Y_t)+\sigma\, dB^C_t.
\end{equation}

It is important to discuss in more detail the trajectory fluctuations that are here being accounted for. If we choose to consider the fluctuations belonging to the normal frame $N_{\gamma_0}M$ over any $V^n$, the dimensionality of the contribution of the fluctuations to the final propagator will be equal to $n-1$. This can be seen by computing (4.5) for a free particle, which yields $(\frac{m}{2\pi i})^{\frac{n-1}{2}}$ rather than the standard n-dimensional result from the traditional path integral. Luckily, this issue can be solved by taking into account longitudinal fluctuations. 
\vspace{1em}

So far, we have been considering the fluctuations $X_t$ as being defined over the normal bundle, which would yield an orthogonal decomposition roughly of the form $X(t)=\gamma_0(t)+p_N(t):p_N \in N_{\gamma_0}$. Now, if we consider longitudinal fluctuations $p_L(t): p_L(0)=p_L(1)=0,$ we could correct our fluctuation to be $X_t=\gamma_0(t)+p_N(t)+p_L(t)\dot{\gamma}_0(t)$. Under infinitesimal reparameterization on the Brownian bridge $t \rightarrow t+\epsilon(t):\epsilon(0)=\epsilon(1)=0,$, we obtain that $\delta X(t)=\epsilon(t)\dot{\gamma}(t) \space \forall \gamma \in \text{Hom} \gamma_0$, wherefrom longitudinal variations are gauges. Therefore, we find ourselves in need of removing the zero mode. This may be done by introducing a gauge condition 
\begin{equation}
G(p_L)=\int_0^1 p_L(t)dt=0,
\end{equation}
We may then introduce the FP identity, by which we get (3.14), removing the zero mode. 
\begin{equation}
    \int \delta(G(p_L+\epsilon))\det (\frac{\delta G}{\delta \epsilon}) \mathcal{D}\epsilon=1
\end{equation}

\section{The Gaussian Expectation as a Stochastic Path Integral}
\vspace{1em}

\subsection{Formulating the Measure and the Stochastic Path Integral}
\vspace{1em}

The stochastic dynamics of individual fibers allow us to formulate a probability law for the evolution of the system. Let $\Gamma$ be the Gaussian law for the covariant bundle Brownian motion on $C([0,T], \gamma^*TM)$ where $T$ is the time it takes for the system to reach B from A along $\gamma$. When we generalize over all paths in $\mathcal{T}_\eta,$ we will cease to view each path as a physical trajectory and view them as curves, which will be parametrized by $[0,1]$. In the parallel frame, we have the law given by (3.8), with the stochastic field $b_t=-\frac{1}{\sigma^2} \nabla U_t(x)$ and $U_t(x)=\frac{1}{2}\chi(t)|x|^2+E_c(|x|) : E_c:[0,\eta ) \rightarrow \mathbb{R}$ being a bounded $C^1$ barrier that is asymptotic at the boundary of the tube. Assuming that $E_c$ is locally Lipshitz on $\text{int} \mathcal{T}_\eta,$ $\exists X_t \in \{ |x| < \eta \},$ with some $X_0$ being a unique solution to the SDE. Let there be a Gaussian measure path law $\mu $ for all $X$. Then it follows from the Girsanov theorem that 
\begin{equation}
    \frac{d \mu}{d \Gamma}(X)=\exp (\frac{1}{\sigma} \int_0^T \langle b_t(X_t), dB_t \rangle -\frac{1}{2}\int^T_0|b_t(X_t)|^2dt)
\end{equation}
On $TM$ with $\pi : E \rightarrow [0,1], \quad E= \cup_{t\in[0,1]}T_{\gamma(t)}M, \quad \exists \mu^{TM}: \forall A$---measurable, $\mu^{TM}=\int_0^1\mu_t(A)\bar{\rho}(dt), $ where $\bar{\rho}=\pi_*\mu^{TM}$ a conditional fiber and $ \quad \mu^{TM}$---a $\sigma-$ additive probability on $E$ and $A$ is Borel, belonging to the Borel $\sigma-$algebra $\mathcal{B}([0,1],M)$. For each $t\in[0,1]$, let
$\mu_t := (\mathrm{ev}_t)_*\mu,\qquad 
\mu_t(A) = \mu\big(\{\omega\in\mathcal C:\ X_t(\omega)\in A\}\big),$ 
which define the $t$---marginal on the fibre $T_{\gamma(t)}M$. Let $E := \bigsqcup_{t\in[0,1]} T_{\gamma(t)}M$ be the bundle with projection $\pi:E\to[0,1]$, and let $\bar\rho$ be a probability measure on $[0,1]$, which is the bundle we have defined for the flux tubes. The Gaussian measure $\Gamma$ splits into a longitudal and normal measures, such that $\Gamma_N \in H^1_0([0,1], \mathbb{R}^{n-1})$ and $\Gamma_L \in H^1_0([0,1])$, with $\Gamma =\Gamma_N\otimes\Gamma_L$. 
\vspace{1em}

Then for every bounded Borel function $F:E\to\mathbb R$,
\begin{equation}
  \int_E F(x,t)\,\mu^{TM}(dx,dt)
  \;=\;
  \int_0^1 \mathbb E_\mu\big[F(X_t,t)\big]\,\bar\rho(dt).
\end{equation}
This identity may be seen through application to the indicator function $\mathbf{1}_A$, for which we have
\begin{equation}
    \int_E \mathbf 1_A(x,t)\,\mu^{TM}(dx,dt)
= \mu^{TM}(A)
= \int_0^1 \mu_t(A_t)\,\bar\rho(dt).
\end{equation}

For the indicator function, we have the identity $\mathbf 1_A(X_t,t)= \mathbf 1_{\{(X_t,t)\in A\}}= \mathbf 1_{A_t}(X_t),$ so by the definition of $\mu_t$ as the law of $X_t$, $\mathbb E_\mu[\mathbf 1_A(X_t,t)]
= \mathbb E_\mu[\mathbf 1_{A_t}(X_t)]
= \mu_t(A_t).$

\begin{equation}
\therefore \int_0^1 \mathbb E_\mu[\mathbf 1_A(X_t,t)]\,\bar\rho(dt)
= \int_0^1 \mu_t(A_t)\,\bar\rho(dt)
= \mu^{TM}(A)
= \int_E \mathbf 1_A(x,t)\,\mu^{TM}(dx,dt)
\end{equation}

This also holds true for all simple functions $F(x,t)=\sum_{k=1}^n a_k\,\mathbf 1_{A_k}(x,t),$ and similarly for all general bounded Borel functions. Under these conditions, we may formulate the Path Integral as the expectation value over $\mu$;
\vspace{1em}

\begin{definition}
$\forall O(X)-$observables, the Stochastic Path Integral is given by the expectation value 
\begin{equation}
    I_\mu(O)=\mathbb{E_\mu}[O(X)\exp (\frac{i}{\hbar}S(X))]=\mathbb{E}_\Gamma[O(X) \frac{d\mu}{d\Gamma}\exp (\frac{i}{\hbar}S(X))]
\end{equation}
\end{definition}

 The rest of the current section will be dedicated to establishing the properties of the Stochastic Path integral. We would also like to underline that the phase factor $\exp(\frac{i}{\hbar}S)$ is posited. While it is possible to derive it from principles of statistical mechanics in application to our setup, this will not be performed here. 
    
\begin{theorem}[Fubini theorem]
Let $(\Omega,I,\mu)$ be the probability space for the stochastic process in $\mathcal{T}_\eta$ with the diffusion $X=(X_t)_{t\in[0,1]}$ and values in the bundle $E=\bigsqcup_{t\in[0,1]}T_{\gamma(t)}M$. Then, for every bounded Borel function $F:E\to\mathbb{R}$,
\begin{equation}
\label{eq:SFPI-Fubini}
  \int_E F(x,t)\,\mu^{TM}(dx,dt)
  \;=\;
  \int_0^1 \mathbb E_\mu\big[F(X_t,t)\big]\,\bar\rho(dt).
\end{equation}
\end{theorem}

\textbf{Proof:} Consider the product probability space $(\Omega\times[0,1],\ I\otimes\mathcal B([0,1]),\ P),
\qquad P:=\mu\otimes\bar\rho.$ Then, the associated measurable map $\Phi:\Omega\times[0,1]\to E,\qquad
\Phi(\omega,t):=\bigl(X_t(\omega),t\bigr).$ is such that for any Borel measurable set $A\subset E$---by definition of $\mu^{TM}$ and $\mu_t$---we have $\mu^{TM}(A)
= \int_0^1 \mu_t(A_t)\,\bar\rho(dt)
= \int_0^1 \mu\big(\{\omega:\ X_t(\omega)\in A_t\}\big)\,\bar\rho(dt)$ with $\{\omega:\ X_t(\omega)\in A_t\}
= \{\omega:\ (X_t(\omega),t)\in A\}
= \{\omega:\ \Phi(\omega,t)\in A\}$, and hence
\begin{equation}
  \mu^{TM}(A)
= \int_0^1 \mu\big(\{\omega:\ \Phi(\omega,t)\in A\}\big)\,\bar\rho(dt)
= \int_{\Omega\times[0,1]} \mathbf 1_{\Phi^{-1}(A)}(\omega,t)\,P(d\omega,dt)
= P\big(\Phi^{-1}(A)\big).  
\end{equation}
Wherefrom, we see that $\mu^{TM} = \Phi_*P,$, implying that $\mu^{TM}$ is the pushforward of $P=\mu\otimes\bar\rho$ under $\Phi$. In the integral representation, let $F:E\to\mathbb{R}$ be bounded and Borel measurable. Since $\mu^{TM}=\Phi_*P$, we have by the definition of pushforward measure
\begin{equation}
    \int_E F(x,t)\,\mu^{TM}(dx,dt)
= \int_{\Omega\times[0,1]} F\bigl(\Phi(\omega,t)\bigr)\,P(d\omega,dt)
= \int_{\Omega\times[0,1]} F\bigl(X_t(\omega),t\bigr)\,\mu(d\omega)\,\bar\rho(dt).
\end{equation}
By the usual Fubini theorem for the product measure $\mu\otimes\bar\rho$, this equals 
\begin{equation}
    \int_0^1 \left(\int_\Omega F\bigl(X_t(\omega),t\bigr)\,\mu(d\omega)\right)\bar\rho(dt)
= \int_0^1 \mathbb E_\mu\big[F(X_t,t)\big]\,\bar\rho(dt),
\end{equation}
which is precisely (4.6). \qed
\vspace{1em}

\subsection{Semigroup structure and Riemann Product}
\vspace{1em}

In the case of the Feyman-Kac theorem and all of the existing rigorous formulations of the Path Integral, the formulae define Markov semigroups via their diffusion generators. In the bundle framework, the bundle law captures all the marginals defined for the same $t \in [0,1]$, and hence to represent the amplitudes over a multitude of points, we may define a transition kernel $K(Q_x,Q_y)$ between two configurations, which given the stochasticity of the diffusion process in $\mathcal{T}_\eta$ is a Markov kernel. Alongside the kernel, we may also define a time-dependent generator of $X$ in the parallel frame of the coordinates for the tube with $G_tf=\langle b_t, \nabla f \rangle +\frac{\sigma^2}{2}\nabla^2f$. 
\vspace{1em}

Let $V$ be the local density defining a twisted Feynman-Kac operator \begin{equation}
    u(s,x)=T_{s,t}f(x)=\mathbb{E}_x[\exp(\frac{i}{\hbar}\int_s^tV(u,X_u)du)f(X_t)|X_s=x]
\end{equation}
By the Feynman-Kac theorem, we have that the generators obey the following differential equation
\begin{equation}
    \frac{\partial u}{\partial s}+G_su=\frac{i}{\hbar}V(s,x)u
\end{equation}
In our case, we may posit a generator incorporating the stochastic coefficients of (3.8)-(3.12) of form $\mathcal{L}f=\partial_tf+G_tf$. It is possible to verify that this diffusion generator also forms a semigroup, and is a generator for a space-time Markov process. Below, we will seize to write $\hbar$ for the sake of practicality

\begin{theorem}
Let $p$ be the partition of $[0,1]$ and $K_{t_i,t_{i+1}}$---Markov kernel of X. Then, the Riemann product for $I_\mu$ is 
\begin{equation}
    I_p=\int\mu_{t_0}(dx_0)\int... \int \prod^{n-1}_{i=0} \exp(i\mathcal{L}(t_i,x_i)\Delta t_i)K_{t_i,t_{i+1}}(x_i,dx_{i+1})\mu_{t_i}(dx_i)
\end{equation}
such that given $|p|\rightarrow0$, $I_p\rightarrow I_\mu$. 
\end{theorem}
\vspace{1em}

\textbf{Proof:} From the definition of the Markov kernels $K_{t_i,t_{i+1}}$ and the Markov property, we may express $I_p$ as
$I_p=\mathbb{E}_\mu\big[e^{\,iS_p}\big],$ so $|I_\mu-I_p|$ follows
\begin{equation}
\label{eq:diff-Imu-Ip}
|I_\mu-I_p|
=\big|\mathbb{E}_\mu\big[e^{\,iS}-e^{\,iS_p}\big]\big|
\le \mathbb{E}_\mu\big[\,|e^{\,iS}-e^{\,iS_p}|\,\big].
\end{equation}
$\forall u,v, \quad |e^{iu}-e^{iv}|
\le |u-v|,$ and so applying this with $u=S(\omega)$, $v=S_p(\omega)$, we get $|e^{\,iS}-e^{\,iS_p}|
\le |S-S_p|.$ It then follows that $|I_\mu-I_p|
\le \mathbb{E}_\mu\big[\,|S-S_p|\,\big].$ By defining $S=\int_0^1\mathcal{L}(t,X_t)\,dt,\quad
S_p=\sum_{i=0}^{n-1}\mathcal{L}(t_i,X_{t_i})\,\Delta t_i,$
we get that for each fixed stochastic path $\omega \in \Omega$ we have the Riemann sum approximation
\begin{equation}
    S(\omega)=\int_0^1\mathcal{L}\big(t,X_t(\omega)\big)\,dt,
\qquad
S_p(\omega)=\sum_{i=0}^{n-1}\mathcal{L}\big(t_i,X_{t_i}(\omega)\big)\,\Delta t_i.
\end{equation}
Since $X_t(\omega)\in\mathcal{T}_\eta$ and $\mathcal{L}$ is continuous on the compact set
$[0,1]\times\mathcal{T}_\eta$, the map $t\rightarrow\mathcal{L}(t,X_t(\omega))$ is continuous, bounded on $[0,1]$. Therefore, $|S(\omega)-S_p(\omega)|\xrightarrow[|p|\to0]{}0,
\quad \forall\,\omega\in\Omega,$ and $\mathcal{L}(t,X_t(\omega))|\le C,\quad\forall t,\omega,$
for some constant $C<\infty$, and hence $|S(\omega)|\le C,|S_p(\omega)|\le C \implies |S(\omega)-S_p(\omega)|\le 2C,\quad\forall\,\omega.$ Thus $|S-S_p|\to0$ pointwise $\mu$–a.s.\ as $|p|\to0$, and $|S-S_p|\le 2C$ is integrable. By DCT, $\lim_{|p|\to0}\mathbb{E}_\mu\big[\,|S-S_p|\,\big]=0.$ This yields $\lim_{|p|\to0}|I_\mu-I_p|
\le \lim_{|p|\to0}\mathbb{E}_\mu\big[\,|S-S_p|\,\big]
=0,$
which implies $I_p\to I_\mu$ as $|p|\to0$.\qed 
\vspace{1em}

It is important to underline that $\mathcal{L}$ above is not the Lagrangian, but the semigroup generator mentioned previously. The generators do act as analogues of the Lagrangians and yield a short-time propagator approximation of form $\exp(i\mathcal{L}\Delta t),$ whose product yields an exponential function with the integral (4.14) inside. 
\vspace{1em}

\section{Equivalence of the Stochastic Path Integral and Euclidean Theory}
\vspace{1em}

The present section will be dedicated to establishing and proving equivalence between the stochastic framework and the existing rigorous formulations of the path integral realized under Wick's theorem, which is valid for the Stochastic Path Integral as well. During the writing of this paper, the notation for the tubular bundle changed from $\mathcal{T}_\eta$ to $T_\eta$ in an attempt to simplify the typing process. We make this notation change explicit and the author apologizes for any and all inconvenience caused by this. Since the Feynman-Kac Path Integral also is based on stochastic dynamics, in some cases the path integral from 4.5 will be referred to as the 'tube-confined'.
\vspace{1em}

\subsection{The Stochastic Path Integral and Wick Rotation}

\begin{lemma}
Under the transformation $t\rightarrow it, \quad I_\mu(O)=\mathbb{E_\mu}[O(X)\exp (\frac{i}{\hbar}S(X))] \rightarrow \mathbb{E_\mu}[O(X)\exp (-\frac{1}{\hbar}S(X))]$. The latter expression exists and converges.
\end{lemma}

\textbf{Proof.} Since $O(X)$ is bounded with $(e^{-S_E},e^{\,iS_E})=\beta$ having a modulus $|\beta|\leq1$, both $I_\mu$ and $I_E=I_\mu^{t\rightarrow it}$ exist as (finite) expectations. For the Euclidean case,
\begin{equation}
|I_E-I_{E,p}|
=\big|\mathbb{E}_\mu\big[O(X)\big(e^{-S_E}-e^{-S_{E,p}}\big)\big]\big|
\le \|O\|_\infty\,\mathbb{E}_\mu\big[\,|e^{-S_E}-e^{-S_{E,p}}|\,\big].
\end{equation}
Since $\forall u,v \in \mathbb{R}, \quad |e^{-u}-e^{-v}|\le|u-v|, \quad |I_E-I_{E,p}|
\le \|O\|_\infty\,\mathbb{E}_\mu\big[\,|S_E-S_{E,p}|\,\big].$
As in Theorem 4., continuity and boundedness of $\mathcal{L}_E$ on the compact tube imply $|S_E(X)-S_{E,p}(X)|\xrightarrow[|p|\to0]{}0\quad\text{for all }X,
\qquad
|S_E(X)-S_{E,p}(X)|\le 2C$ for some constant $C$, so again by DCT,
$\mathbb{E}_\mu[|S_E-S_{E,p}|]\to0$ meaning that $I_{E,p}\to I_E$.  Thus we get $|I_\mu-I_{\mu,p}|
\le \|O\|_\infty\,\mathbb{E}_\mu\big[\,|S_E-S_{E,p}|\,\big]\xrightarrow[|p|\to0]{}0,$ so $I_{\mu,p}\to I_\mu$. We get that the Wick–rotated stochastic path integral exists and equal to the limit of the the Stochastic Path Integral in representation (4.12). \qed
\vspace{1em}

Under this theorem, arguments that relate the Stochastic Path Integral to the rigorously defined Euclidean analogue gain legitimacy. As such, it is possible to relate the current framework to existing work on a rigorous formulation of the Path Integral, notably the Feynman-Kac theorem. 
\vspace{1em}

The transition to complex time also however modifies our stochastic differential equations in a way that leads $(dW_t)^2=idt$. Due to this failure, we obtain a necessity to alter our stochastic dynamics in a manner that remains meaningful relative to the noise, which may be found via the complexification of diffusion generators. This may be done by the passage from a real Kato potential $V$ to a complex $iV$. This gives us a diffusion generator of form $\mathcal{L}_\mathbb{C}=\partial_t+G_t+iV$.
\vspace{1em}

\subsection{Stochastic Path Integral and the Feynman-Kac theorem}
\vspace{1em}

Let $V:[0,T] \times D \rightarrow \mathbb{R}$ be bounded and measurable for $D \subset \mathbb{R}^n$. For any $f:D \rightarrow \mathbb{R}, \quad \forall x \in D, \space f(x)\leq C,$ we may define the expectation value $u(s,x)$
\begin{equation}
    u(s,x)=\mathbb{E}_{s,x} [\exp(-\int_s^TV(r,X_r)dr)f(X_t)]
\end{equation}
The Feynman-Kac Theorem states provides a rigorous general grounding for the Path integral, stating that $u(s,x)$ is a unique, bounded solution of the differential equation $\partial_su(s,x)+L_su(s,x)-V(s,x)=0$ with $u(T,x)=f(x)$ under $s<T, x\in D.$ Here, $L_s$ is a diffusion generator that will not be specified. This theorem was very broadly generalized in \cite{FK}. Under the passage to the $\eta-$Tube representation and the complexified Stochastic Path Integral framework, we expect to obtain the following 2 theorems;
\vspace{1em}

\begin{theorem} Let $X=(X_t)_{t\in[0,T]}$ be the diffusion of $T_\eta(\gamma_0)$ and let $V_E:[0,T]\times T_\eta(\gamma_0)\to\mathbb{R}$ be the Euclidean density of the
action such that $\sup_{(t,x)\in[0,T]\times T_\eta(\gamma_0)} |V_E(t,x)| =: C < \infty.$ Let $f:T_\eta(\gamma_0)\to\mathbb{C}$ be bounded Borel, and for $\theta\in\mathbb{C}$ and
$x\in T_\eta(\gamma_0)$. Then, the $\theta-$modified version of (5.2) is 
\begin{equation}\label{eq:def-u-theta}
  u_\theta(0,x)
  := \mathbb{E}_{0,x}\!\left[
    f(X_T)\exp\!\Big(-\theta\int_0^T V_E(t,X_t)\,dt\Big)
  \right].
\end{equation}
For every $x\in T_\eta(\gamma_0)$, the expectation in exists and is finite $\forall \theta\in\mathbb{C}$ and the map $\theta\mapsto u_\theta(0,x)$ extends to an entire function on $\mathbb{C}$.
\end{theorem}
\vspace{1em}

\textbf{Proof:} Let us fix some concrete $x\in T_\eta(\gamma_0)$. Then, we may write the null $x$ expectation $\mathbb{E}_{0,x}$ by setting 
\begin{equation}
    A := \int_0^T V_E(t,X_t)\,dt.
\end{equation}
wherefrom by the boundedness assumption on $V_E$ we have a pathwise-defined inequality,
\begin{equation}
  |A|\le \int_0^T |V_E(t,X_t)|\,dt \le \int_0^T C\,dt = CT.
\end{equation}
Furthermore, since $f$ is bounded, $\exists M_f<\infty$ with $|f(X_T)|\le M_f.$
Now, we shall prove the existence and finiteness of the expectation  $\forall \theta\in\mathbb{C}$. For some arbitrary $\theta$, $\left|\exp(-\theta A)\right|
= \exp\big(\mathrm{Re}(-\theta A)\big)
\le \exp\big(|\theta||A|\big)
{\le} \exp\big(|\theta|CT\big).$ Hence, $\big| f(X_T)\exp(-\theta A)\big|
\le M_f \quad \exp(|\theta|CT)
= M(\theta),$ where $M(\theta)<\infty$ depends only on $\theta$, $C$, $T$ and $\|f\|_\infty$. In particular, the integrand is bounded and therefore integrable with respect to some
probability measure $\mathbb{P}_{0,x}$. Thus, the expectation for $\theta$ exists and is finite $\forall\theta\in\mathbb{C}$. Finally, for each $\omega$ in the underlying probability space, the map $\theta \longmapsto f(X_T(\omega))\,e^{-\theta A(\omega)} $ is a holomorphic function of $\theta$ on the entirety of the complex plane. Such a function is referred to in literature as entire. To show that $\theta\mapsto u_\theta(0,x)$ is entire, we may represent it with a power series which ought to have a radius of convergence at infinty. For any $\theta\in\mathbb{C}$ and any $\omega$,
\begin{equation}
    e^{-\theta A(\omega)}
= \sum_{n=0}^\infty \frac{(-\theta A(\omega))^n}{n!}
\end{equation}
This series is absolutely convergent, as may be verified below; 
\begin{equation}
    \sum_{n=0}^\infty \left|\frac{(-\theta A(\omega))^n}{n!}\right|
\le \sum_{n=0}^\infty \frac{(|\theta||A(\omega)|)^n}{n!}
{\le}
\sum_{n=0}^\infty \frac{(|\theta|CT)^n}{n!}
= e^{|\theta|CT}.
\end{equation}
Multiplying by $f(X_T(\omega))$ and leveraging the established bound $|f(X_T(\omega))|\le M_f$, we obtain the pointwise bound
\begin{equation}
\sum_{n=0}^\infty \left|
  f(X_T(\omega))\frac{(-\theta A(\omega))^n}{n!}
\right|
\le M_f e^{|\theta|CT}.
\end{equation}
For any fixed radius $R$ that majorates $\theta,$
\begin{equation}
\sum_{n=0}^\infty \left|
  f(X_T(\omega))\frac{(-\theta A(\omega))^n}{n!}
\right|
\le M_f e^{RCT}
\end{equation}
Therefore, by dominated convergence we may interchange expectation and summation, obtaining the following;
\begin{equation}
u_\theta(0,x)
= \mathbb{E}\big[f(X_T)e^{-\theta A}\big]
= \mathbb{E}\left[\sum_{n=0}^\infty f(X_T)\frac{(-\theta A)^n}{n!}\right]
= \sum_{n=0}^\infty \mathbb{E}\big[f(X_T)A^n\big]\frac{(-\theta)^n}{n!},
\end{equation}
Let $a_n:=\mathbb{E}[f(X_T)A^n]$. From established results, it follows that $|a_n|
\le \mathbb{E}\big[|f(X_T)||A|^n\big]
\le M_f (CT)^n.$ Hence
\begin{equation}
    \limsup_{n\to\infty} \sqrt[n]{\frac{|a_n|}{n!}}
\le \limsup_{n\to\infty} \sqrt[n]{\frac{M_f(CT)^n}{n!}} = 0,
\end{equation}
so the power series for our entire function has infinite radius of convergence. Therefore $\theta\mapsto u_\theta(0,x)$ is an entire
function on $\mathbb{C}$. \qed
\vspace{1em}

We may underline here that the Feynman-Kac theorem may be defined on the $\eta$ tube. This is easy to see given that $\forall \theta >0,$ with $0 \leq t \leq T, $ the $\theta$ expectation (5,3) on $(t,x)$ is given by $u_\theta:[0,T] \times T_\eta \rightarrow \mathbb{C}$ which is a unique, bounded classical solution of the equation 
\begin{equation}
    \frac{\partial u_\theta}{\partial t}+G_tu_\theta - \theta V_E u_\theta =0
\end{equation}
Here, $u_\theta(T,x)=f(x):x \in T_\eta$. This presents a first instance of the application of the (real) Feynman-Kac theorem to the tubular neighborhood construction. This may be elaborated on in the following theorem;

\begin{theorem}
Let $f:T_\eta(\gamma_0)\to\mathbb{C}$ be bounded Borel, such that for $\theta\in\mathbb{C}$ and
$x\in T_\eta(\gamma_0)$ the expectation value (5.3) exists. Then, for $\theta=-i$ the expectation exists as an absolutely convergent Lebesgue integral $\forall x\in T_\eta(\gamma_0)$, such that if $V_E=\tfrac{1}{\hbar}L_E$ is the Euclidean Lagrangian density and
$S(X)=\int_0^T L(t,X_t,\dot X_t)\,dt$ is the corresponding action
obtained by Wick rotation, then
\begin{equation}
    \exp\!\Big(i\int_0^T V_E(t,X_t)\,dt\Big)
=\exp\!\Big(\tfrac{i}{\hbar}S(X)\Big)
\end{equation}
with $u_{-i}(0,x)$ coinciding with the tube-confined stochastic Feynman path integral with observable $O(X)=f(X_T)$.
\end{theorem}
\vspace{1em}

\textbf{Proof:} We have the same set up as in the proof of Theorem 6. By the boundedness of $V_E$, we once again recover the same pathwise-defined inequality
\begin{equation}
  |A|\le \int_0^T |V_E(t,X_t)|\,dt \le \int_0^T C\,dt = CT.
\end{equation}
Since $f$ is bounded, $\exists M_f<\infty$ with $|f(X_T)|\le M_f$. For $\theta=-i$ we obtain that $\Big|f(X_T)e^{iA}\Big|
= |f(X_T)|\,\big|e^{iA}\big|
= |f(X_T)| \le M_f,$
since $|e^{iA}|=1 \quad \forall A \in \mathbb{R}$. In particular, the integrand in for the expectation under $\theta=-i$ is bounded and hence integrable with respect to some probability
measure $\mathbb{P}_{0,x}$, wherefrom the expectation exists
and is finite, while the stochastic Feynman weight $\exp({i\int_0^T V_E})$ 
in the integrand admits a finite-dimensional Lebesgue integral without necessitating an oscillatory regulation.
Assuming that $V_E=L_E$ with $L_E$ the Euclidean
Lagrangian density, we get that under
Wick rotation (as fixed in the main text), the Euclidean action
$\int_0^T L_E\,dt$ corresponds to the  action $S(X)=\int_0^T L(t,X_t,\dot X_t)\,dt$,
and we obtain
\begin{equation}
    \int_0^T V_E(t,X_t)\,dt
= \int_0^T L_E(t,X_t,\dot X_t)\,dt 
= S(X) \quad \exp\!\Big(i\int_0^T V_E(t,X_t)\,dt\Big)
= \exp\!\Big(\tfrac{i}{\hbar}S(X)\Big).
\end{equation}

By the substitution of the expectation for $\theta=-i$, we get
\begin{equation}
    u_{-i}(0,x)
= \mathbb{E}_{0,x}\!\left[
  f(X_T)\exp\!\Big(iS(X)\Big)
\right],
\end{equation}
which is precisely the tube-confined stochastic Feynman path integral with observable
$O(X)=f(X_T)$ under $\hbar=1$. \qed

\section{Derivation of the Free Particle Propagator}

From a theoretical standpoint, the framework constructed above provides a robust measure and an object for path space compactification for a free particle with a quadratic Lagrangian traveling between two points with finite energy uncertainty and time of travel. It is not necessary to verify whether the framework holds by computing the propagator for a free particle with mass $m$ and a base-trajectory $\gamma_0(t)=A+\frac{t}{T}(B-A):T<\infty, t\in [0,1]$ in $\mathbb{R}^n$. 
\vspace{1em}

The trajectories in the tube are given by $X_t=\gamma_0(t)+Y_t:Y_0=Y_T=0$, with the action for the system being defined by the integral
\begin{equation}
    S(\gamma)=\frac{m}{2}\int_0^T \dot{\gamma}(t)^2dt
\end{equation}
From $\dot{\gamma}_0(t)=\frac{|B-A|}{T}(1-t),$ we find that $S(\gamma_0)=\frac{m}{2T}|B-A|^2\implies S(\gamma)=\frac{m}{2T}|B-A|^2+\frac{m}{2T}\int_0^T|\dot{Y}_t|^2dt$. From this, the propagator for the function is simply 
\begin{equation}
    K(B,T;A,0)=\exp(\frac{im}{2T\hbar}|B-A|^2) \times \mathbb{E}_\mu[\exp(\frac{im}{2 \hbar}\int_0^T|\dot{Y}_t|^2dt)]
\end{equation}
We have previously seen that by Girsanov, we may define a Gaussian measure $\Gamma$ on $H^1_0([0,T];\mathbb R^n)$ satisfying (4.1). Then, under $\Gamma$, $Y\in H^1_0:\int_0^T|\dot{Y}_t|^2dt <\infty \implies \exp(\frac{im}{2 \hbar}\int_0^T|\dot{Y}_t|^2dt)=F(Y):|F(Y)|=1, \quad F\in L^1(\Gamma) \cap L^\infty (\Gamma).$ We may remove the prefactor of $\frac{im}{2 \hbar}$ and just view this expression through the lens of a random variable under $\Gamma$, 
\begin{equation}
   Q(Y):=\int_0^T |\dot Y_t|^2\,dt 
\end{equation}

This is a quadratic form on the real Hilbert space $H:=H^1_0([0,T];\mathbb R^n)$ with a corresponding inner product $\langle Y,Z\rangle_H
:=
\int_0^T \dot Y_t\cdot \dot Z_t\,dt.$
Meanwhile, the Gaussian measure $\Gamma$ has covariance operator $C = (-\Delta_D)^{-1},$
where $-\Delta_D$ is the Dirichlet Laplacian on $[0,T]$. Then, let there be some family of functions $\{ e_k(t) \}$ with 

\begin{equation}
    e_k(t)
=
\sqrt{\frac{2}{T}}
\sin\!\left(\frac{k\pi t}{T}\right),
\qquad
\lambda_k=\left(\frac{k\pi}{T}\right)^2.
\end{equation}

Then, under $\Gamma$,
\begin{equation}
    Y(t)
=
\sum_{k=1}^\infty \xi_k e_k(t),
\end{equation}

with $\xi_k \sim \mathcal N\!\left(0,\lambda_k^{-1}\right)$ defined component-wise in $\mathbb R^n$. We may then compute $Q(Y)
=
\sum_{k=1}^\infty \lambda_k\,|\xi_k|^2,$ obtaining that 
\begin{equation}
    \mathbb E_\Gamma\!\left[
e^{\frac{i m}{2\hbar}Q(Y)}
\right]
=
\prod_{k=1}^\infty
\mathbb E\!\left[
\exp\!\left(
\frac{i m}{2\hbar}\lambda_k|\xi_k|^2
\right)
\right]^{n}.
\end{equation}

For a Gaussian random variable over $\mathbb{R}$, $\xi\sim\mathcal N(0,\lambda^{-1})$, $\mathbb E\!\left[e^{i a \lambda \xi^2}\right]
=
(1-2ia)^{-1/2},
\qquad a\in\mathbb R.$
Hence,
\begin{equation}
    \mathbb E\!\left[
e^{\frac{i m}{2\hbar}\lambda_k\xi_k^2}
\right]
=
\left(1-\frac{i m}{\hbar}\right)^{-1/2}.
\end{equation}

Taking the infinite product yields
\begin{equation}
    \mathbb E_\Gamma\!\left[
e^{\frac{i m}{2\hbar}Q(Y)}
\right]
=
\prod_{k=1}^\infty
\left(
1-\frac{i m}{\hbar}
\right)^{-n/2}.
\end{equation}

This product is understood in the zeta–regularized sense, giving the determinant representation
\begin{equation}
    \mathbb E_\Gamma\!\left[
e^{\frac{i m}{2\hbar}Q(Y)}
\right]
=
\det\nolimits^{-n/2}
\!\left(
I-\frac{i m}{\hbar}(-\Delta_D)^{-1}
\right).
\end{equation}

By leveraging the classical identity $\prod_{k=1}^\infty
\left(
1+\frac{a^2}{k^2\pi^2}
\right)
=
\frac{\sinh a}{a},$ 
together with analytic continuation $a\mapsto e^{-i\pi/2}a$, we obtain
\begin{equation}
    \mathbb E_\Gamma\!\left[
\exp\!\left(
\frac{i m}{2\hbar}\int_0^T |\dot Y_t|^2\,dt
\right)
\right]
=
\left(\frac{m}{2\pi i\hbar T}\right)^{\frac{n}{2}}
\end{equation}

This last expression is a genuine $L^1(\Gamma)$ expectation, not a distributional limit, and it is evaluated entirely in real (Lorentzian) time. Therefore, under Gaussian time, we find a result equal to the results of the traditional path integral for the free particle.

\section{Conclusion and Prospective directions}
\vspace{1em}

We have formulated a restricted path space that encompasses all physically meaningful trajectories which may be endowed with a Hilbert space bundle structure and a stochastic process, yielding a path integral that does not face the issues of divergence plaguing the heuristic path integral with an undefined measure $\mathcal{D}x(t)$. The framework as detailed in the current paper, however, is still very raw and bound to the preliminary conditions detailed in the introduction.
\vspace{1em}

As mentioned in the introduction, one of the most significant problems of the path integral is the issue of convergence, which depends on rapid oscillations at infinity. This issue underpins attempts of correcting the path integral through cycles in the complex plane, realized in works by Witten et al, which may be seen in \cite{WIT}. In the context of the Stochastic Path Integral as given in (4.5), the problem of convergence persists, yet due to the bound in the difference of action and the existence of a maximal value for $S$, the phase factor is at all times bounded, and hence the issue of divergence due to that is removed similarly to how the destructive interference acts in the case of the Feynman formulation.
\vspace{1em}

Upon considering a space of all possible tubes that correspond each to neighborhoods of trajectories determined by the energies in the spectrum of the particle, we may still however run into problems of convergence. This may be solved by the means of some 'summation over topologies', where the topologies are the bundles of each classical trajectory tube. 
\vspace{1em}

Additionally, the path integral detailed here may be extended to gauge theories through additional phase factors in terms of stochastic integrals such as Stratonovich integrals. Nevertheless, the present framework lays ground for a foundation of what may be a new approach to viewing the path integral.
\vspace{1em}

\appendix
\section{Appendix: On Relevant Trajectories}
\vspace{1em}

In order to define what trajectories we can account for under an uncertainty in energy, we ought to define a mathematical criterion that could quantitatively assess whether a trajectory is allowable--as in the energy uncertainty permits a classical particle to deviate from the classical trajectory by this alternative--or whether that trajectory is forbidden. 
\vspace{1em}

Let us denote the energy uncertainty by $\delta E$ and our probe function by $\Phi$. Given that it is irrelevant whether this function is defined with any specific orientation, we could select it as being defined on a contour formed the classical trajectory $\gamma_0$ and any $\gamma$ in its homotopy class. Naturally, under the supposition that the endpoints of the trajectories are unaltered, the two are sufficient to form a closed contour on path space. We will define $\Phi(\gamma, \gamma_0): \forall |\Delta E| > \delta E, \space \Phi \rightarrow \infty $, while for $|\Delta E| \leq \delta E$, $0<\Phi<\infty $. 
\vspace{1em}

A Natural choice for such a function would be a contour integral with the integrand that is a meromorphic function that has poles at $E_0 \pm \delta E.$ For such a function, we can write
\begin{equation}
    f(E)=\frac{1}{(\delta E)^2 -(\Delta E)^2} \implies f(H)=\frac{1}{(\delta E)^2-(H-E_0)^2}
\end{equation}

The corresponding probe function on can be chosen to also depend on the canonical 1-form to be translation-invariant, thus having $\Phi$ vanish for paths that have a purely gauge nature. We obtain
\begin{equation}
    \Phi(\gamma,\gamma_0)=\oint_{\gamma \circ \gamma_0}f(H)p_idq^i
\end{equation}

For any $\gamma:|H(\gamma)-E_0|> \delta E,$ the contour approaches a pole of $f(H)$, wherefrom by the Residue Theorem $\Phi(\gamma,\gamma_0)$ diverges. The class of admissible trajectories is therefore given by the homotopies of $\gamma_0$ for which $\Phi$ is finite. This function is not rigorous, but is useful for heuristic intuition motivating the restriction of the path space.

\end{document}